\documentclass[a4paper,11pt,fceqn]{article}
\usepackage{latexsym}
\usepackage{amsmath,amssymb}
\usepackage{amsthm}
\usepackage{amsfonts}
\usepackage[usenames]{color}
\renewcommand\baselinestretch{1}
\numberwithin{equation}{section}
\newtheorem{definition}{Definition}[section]
\newtheorem{theorem}{Theorem}[section]

\newtheorem{remark}{Remark} [section]

\allowdisplaybreaks[4]

\begin{document}

\title{{\textbf{Three spheres inequalities and  unique continuation for a three-dimensional Lam\'{e} system of
elasticity with C$^1$ coefficients}}{\footnote{This work was
supported in part by NSFC(No. 10671040), FANEDD (No. 200522), NCET
(No. 06-0359).}}}
\author{{HANG YU}{\footnote{School of
Mathematical Sciences, Fudan University, Shanghai 200433, People's
Republic of  China. Email: hangyumath@hotmail.com}}}
\date{}
\maketitle

$${\textbf {ABSTRACT}}$$

In this paper, a quantitative estimate of  unique continuation is
proved for a  three-dimensional Lam\'{e} system  with $C^1$
coefficients in the form of three spheres inequalities. The property
of the non faster than exponential vanishing of nonzero local
solutions is also given as an application of the three spheres inequality.\\

{\section*{{\normalsize \bf 1. Introduction }}}
\setcounter{section}{1} \setcounter{equation}{0}

We study the three spheres inequality for Lam\'e systems of
elasticity with $C^1$ coefficients. It is a quantitative estimate of
weak unique continuation.    Firstly, let us introduce the Lam\'{e}
system. Assume that $\Omega$ is a bounded domain in $\mathbb{R}^3$,
the Lam\'{e} moduli $\mu=\mu(x)$ and $\lambda=\lambda(x)$ are  $
C^1(\overline{\Omega})$ and satisfy the strong ellipticity
conditions
\[\mu\geq\alpha_0>0,\ \ 2\mu+\lambda\geq\beta_0>0,\] for given positive constants $\alpha_0$ and $\beta_0$.
Generally, it can be  assumed that $\Omega$ contains the origin and
$B_R\subset\subset\Omega$ for some $R>0$ where $B_R$ is an open ball
centered at the origin with radius $R$. The Lam\'{e} system is given
by
\begin{equation}
{\rm{div}}{\big(}\mu(\nabla u+(\nabla
u)^\top){\big)}+\nabla{\big(}\lambda{\rm div}\ u{\big)}=0,
\end{equation}
where $u=(u_1,u_2,u_3)^\top$ is the displacement vector and
\[\nabla u=\left(
             \begin{array}{cccc}
               \partial_{x_1}u_1 &  \partial_{x_2}u_1& \partial_{x_3}u_1  \\
             \partial_{x_1}u_2&  \partial_{x_2}u_2 & \partial_{x_3}u_2 \\
               \partial_{x_1}u_3& \partial_{x_2}u_3 & \partial_{x_3}u_3 \\

             \end{array}
           \right)
\] is the gradient matrix of $u$.

The history of three spheres inequalities is closely related to the
strong continuation uniqueness principle (SUCP). The three spheres
type inequality for scalar elliptic equations is a rather classical
result, dated back to Landis \cite{14}, who generalized the famous
Hadamard's three circles theorem to solutions of elliptic equations.
The three spheres inequality with integral norms for scalar elliptic
equations was originally introduced by Garofalo and Lin \cite{10},
\cite{11} and later developed by Kukavica \cite{13}. Their proof is
based on the monotonicity property of the frequency function
\cite{10},  \cite{13}. Recently, Alessandrini and Morassi \cite{1}
have obtained the three spheres inequality for the isotropic
elasticity system. Their method can be stated as follows. Set a
$(n+1)$-vector valued function
\[U=\left(
                                                        \begin{array}{c}
                                                          u \\
                                                          {\rm{div}}\ u \\
                                                        \end{array}
                                                      \right)\]
where $u$ satisfies  $(1.1)$. The system (1.1) with $C^{1,1}$
coefficients can be reduced to a weakly coupled elliptic system with
Laplacian principal part  \cite{1}, \cite{3}, \cite{8}:
   \[-\Delta U+B(\nabla U)+V U=0,\ \ {\rm{in}}\ D,\]
where $D$ is a bounded domain in $\mathbb{R}^n$.  The coefficient
tensors $B$ and $V$ uniquely depending on $\mu$ and $\lambda$ are
bounded measurable. Then, using Rellich's identity \cite{17}, one
can prove that its corresponding frequency function is monotonous.

Unfortunately, the method mentioned above  can not be used when the
Lam\'{e} coefficients are $C^1$, because  one cannot apply
divergence to $(1.1)$ to diagonalize the system in this case. Unlike
the approach used in \cite{1},  C. -L. Lin et al. derived the three
spheres inequality for a 2-dimensional elliptic system with
$W^{1,\infty}$ coefficients by another type of reduction of the
Lam\'{e} system and  Carleman estimates \cite{15}. It is carried out
by using an auxiliary function $\partial_{x_1}u+T
\partial_{x_2}u$ with an appropriate matrix $T$. The key point is
that their new system contains only first order derivative of the
Lam\'{e} coefficients. However such a reduction may not be applied
to higher dimensions.

Eller  proposed another way to reduce the Lam\'{e} system for the
case of three dimensions \cite{7}. Set
\[ A(\partial)(u_1,u_2)=(\nabla\times u_1+\nabla u_2,-\nabla\cdot
u_1),\]
\[A_\alpha(x,\partial)(u_1,u_2)=(\nabla\times u_1+\alpha\nabla u_2,-\nabla\cdot u_1),\]
where $u_1$ is a vector-valued function with three components, $u_2$
and $\alpha$ are scalar-valued functions.  Choosing $\alpha =
(2\mu+\lambda)/\mu$ and $u_2=0$,  we then have
\[{\big(}\mu\Delta u_1+(\lambda+\mu)\nabla{\rm div}\ u_1,0{\big)}=-\mu A_\alpha(x,\partial)A(\partial)(u_1,0).\]
By transforming the Lam\'{e} system into the composition of  two
first order elliptic operators, the Carleman estimate of the
Lam\'{e} operator with $C^1$ coefficients can be given by the
Carleman estimates of $A(\partial)$ and $A_\alpha(x,\partial)$
(Eller \cite{7}). Then the three spheres inequality for three
dimensional Lam\'{e} system with $C^1$ coefficients can be proved
accordingly.

We state one of the main results of the paper  as follows.

\begin{theorem}
Let  $\Omega$ be a bounded domain in $\mathbb{R}^3$ and the Lam\'{e}
moduli $\mu,\lambda\in C^1(\Omega)$ satisfy the strong elliptic
conditions. For any $R_1,R_2$, $0<R_1<R_2<R$,
\begin{equation}
\int_{B_{R_2}}|u|^2dx\leq
C{\big{(}}\int_{B_{R_1}}|u|^2dx{\big{)}}^{\sigma}{\big{(}}\int_{B_{R}}|u|^2dx{\big{)}}^{1-\sigma}
\end{equation}
holds for $u\in H^2(\Omega)$ being a solution of $(1.1)$, where $C$
and $\sigma\in(0,1)$ are two constants depending on $\frac{R_1}{R}$,
$\frac{R_2}{R}$, $\|\lambda\|_{C^1(\Omega)}$ and
$\|\mu\|_{C^1(\Omega)}$.
\end{theorem}

Quantitative estimates like $(1.2)$ have been shown to be extremely
useful in the treatment of the unique continuation principle and the
inverse boundary value problems \cite{2}, \cite{12}, \cite{13},
\cite{15}. Another result of our paper is related to the  strong
unique continuation. Before stating it we recall a relevant
definition.

\begin{definition} A function $u\in L^2_{loc}(\Omega)$
is said to vanish of infinite order at $x_0\in\Omega$ if  for every
$K\in\mathbb{N}$,
\begin{equation}\int_{|x-x_0|<r}|u|^2dx=O(r^K),\ \ {\rm{as}}\ r\rightarrow0^+.\end{equation}
\end{definition}
\begin{definition}Let $u\in H^2(\Omega)$ be a solution of $(1.1)$. The Lam\'{e} system $(1.1)$ is said to have the strong
unique continuation property  if $u$ satisfies the property that if
there exists a point $x_0\in\Omega$ such that $u$ vanishes of
infinite order at $x_0$, then $u\equiv0$ in $\Omega$.
\end{definition}

The result of the weak unique continuation for the Lam\'{e} system
was first given by Dehman and Robbiano  for $\lambda,\mu\in
C^{\infty}(\mathbb{R}^n)$ \cite{6}.  They proved the Carleman
estimate by pseudodifferential calculus. Then Ang, Ikehata, Trong
and Yamamoto  gave a result for $\lambda\in C^2(\mathbb{R}^n),\
\mu\in C^{3}(\mathbb{R}^n)$ \cite{3};  Weck  proved a result for
$\lambda,\mu\in C^{2}(\mathbb{R}^n)$ \cite{18}, \cite{19}. On the
other hand, the result on the strong unique continuation (SUCP) for
the Lam\'{e} system was first obtained by Alessandrini and Morassi
for $n\geq2$, $\lambda,\mu\in C^{1,1}(\mathbb{R}^n)$ \cite{1}.  Then
Lin and Wang studied the SUCP in the case of $n=2$, $\lambda,\mu\in
W^{1,\infty}(\mathbb{R}^n)$ \cite{16}. Their proof relies on
reducing the Lam\'{e} system to a first order elliptic system and on
some suitable Carleman estimates with polynomial weights.  Recently,
Escauriaza \cite{9} has  proved the SUCP in the case of $n=2$,
$\lambda$ being measurable and $\mu$ being Lipschitz  by a similar
method as the one proposed by Lin and Wang.

In this  paper,  the UCP for the Lam\'{e} system of elasticity will
be proved  for $n=3$ and $\lambda,\mu\in C^1(\Omega)$. The following
theorem is stronger than the weak unique continuation property but a
little weaker than the strong unique continuation property.
\begin{theorem}Assume that   $\Omega$ is a bounded domain in
$\mathbb{R}^3$,
 the Lam\'{e} moduli $\mu,\lambda$ satisfy the strong elliptic
conditions and $u\in H^2(\Omega)$ be a solution to $(1.1)$.

{\rm (i)} Let $\lambda, \mu\in C^1(\Omega)$. If there is a point
$x_0\in\Omega$ and $\varepsilon>0$ such that,
\begin{equation}\int_{|x-x_0|<r}|u|^2dx=O(e^{-{r^{-\varepsilon}}}),\ \ {\rm{as}}\
r\rightarrow0^+.\end{equation} then $u\equiv0,\ {\rm in}\ \Omega.$

{\rm (ii)} Let $\lambda, \mu\in C^2(\Omega)$. If  there is a point
$x_0\in\Omega$ and $\varepsilon>0$ such that,
\begin{equation}\int_{|x-x_0|<r}|\nabla u|^2dx=O(e^{-{r^{-\varepsilon}}}),\ \ {\rm{as}}\
r\rightarrow0^+.\end{equation} then $u\equiv {\rm const.},\ {\rm
in}\ \Omega.$

\end{theorem}
The plan of this paper is  as follows. In Section 2, we will show
the conditionally stability estimate in the Cauchy problem for
$(1.1)$. In Section 3,   the three spheres inequality will be proved
based on the results in Section 2. The unique continuation will be
given in Section 4 as an application of the three spheres
inequality. Throughout the paper, $C$ stands for a generic constant
and its value may vary from line to line.

{\section*{{\normalsize \bf 2. Conditional stability}}}
\setcounter{section}{2}\setcounter{equation}{0}\setcounter{lemma}{0}\setcounter{proposition}{0}\setcounter{theorem}{0}
The Carleman estimate is a powerful technique not only for the
unique continuation, but also for solving the exact controllability,
stability  and the inverse problems.  Carleman estimates are
available for scalar elliptic operators with $C^1$-coefficients
whereas many of the results for elliptic systems require
coefficients with higher regularity.

We consider the equilibrium system

\begin{equation}
\widetilde{L}u(x)=\rm{div}{\big(}\mu(x)(\nabla u(x)+(\nabla
u(x))^\top){\big)}+\nabla{\big(}\lambda(x)\rm{div}\ u(x){\big)}=0.
\end{equation}
Define the Lam\'{e} operator as follows
\begin{equation}
L=\mu(x)\Delta+(\lambda(x)+\mu(x))\nabla\rm{div}.
\end{equation}
Then
\begin{equation}
\widetilde{L}u=Lu+{\big(}\nabla u+(\nabla
u)^\top{\big)}\nabla\mu+(\rm{div}\ u)\nabla\lambda.
\end{equation}

The Carleman estimate of operator $L$ was first given by Dehman and
Robbiano \cite{6} when $\lambda$ and $\mu$ are infinite
differentiable. In this section, we  introduce a Carleman estimate
given by Eller \cite{7} at first, which plays an essential role in
proving the conditional stability of the Cauchy problem for system
$(1.1)$. In his recent work, Eller  proved a Carleman estimate for a
certain first order elliptic system which can be used to prove the
Carleman estimate for the isotropic Lam\'{e} system with
$C^1$-coefficients.\\

\begin{theorem}[Eller]Let $\psi\in C^2(\overline{\Omega})$ have
non-vanishing gradient and set $\varphi=e^{s\psi}-1$ for some $s>0$.
Furthermore, assume that $\mu,\lambda\in C^1(\overline{\Omega})$
satisfy the strong elliptic conditions. Then there exist positive
constants $s_0$ and $C$ such that for $s>s_0$, $\tau>\tau_0(s)$
\begin{equation}
\begin{array}{ll}
\displaystyle\tau^2s^4\int_{\Omega}e^{2s\psi}e^{2\tau\varphi}|u|^2dx+&\displaystyle
s^2\int_{\Omega}e^{2\tau\varphi}|\nabla
u|^2dx\\&+\displaystyle\frac{1}{\tau^2}\int_{\Omega}e^{-2s\psi}e^{2\tau\varphi}|\nabla^2
u|^2dx\leq C\int_{\Omega}e^{2\tau\varphi}|Lu|^2dx
\end{array}
\end{equation}
for $u\in H^2(\Omega)$ with compact support in $\Omega$.
\end{theorem}

The proof of the above theorem can be seen in \cite{7}. Compared to Dehman and Robbiano's
method, Eller's proof is quite simple since no pseudo-differential
calculus is used.

\begin{remark}
From the proof of the above theorem, we know that the constant $C$
depends on $\inf|\mu|$, $\inf|\lambda+\mu|$,
$\|\mu\|_{C^1(\Omega)}$, $\|\lambda\|_{C^1(\Omega)}$ and the weight
function $\psi$'s $C^2$ norm in $\Omega$, but not on $s$ and $\tau$.
\end{remark}

 Now we can get a theorem of the conditional
stability.

\begin{theorem}
Assume that  $\mu,\lambda\in C^1(\Omega)$ satisfy the strong
elliptic conditions and $B_{\theta}\subset B_R$ for some
$\theta\in(0,R)$. Let $\gamma=\partial B_{\theta}$ and
$G=B_R\setminus \overline{B_{\theta}}$. Suppose that $u\in H^2(G)$
solves the Cauchy problem
\begin{equation}\left\{
\begin{array}{lll}
\widetilde{L}u&=&0, \hspace{1.5cm}  {\rm{in}} \hspace{0.3cm}G,\\
\partial^{\alpha}u|_{\gamma}&=&f_{\alpha},\hspace{1cm} |\alpha|\leq1
\end{array}\right.\end{equation}
with $f_{\alpha}\in H^{\frac{3}{2}-|\alpha|}(\gamma)$. Then there
exist a sub-domain $\omega\subset G$ with $\gamma\subset\partial
\omega$ and constants $C>0$, $0<\epsilon<1$ such that
\begin{equation} \|u\|_{L^2(\omega)}\leq C
M_0^{1-\epsilon}\zeta_0^{\epsilon} \end{equation} where
$M_0:=\|u\|_{H^1(G)}$,
$\zeta_0:=\sum_{|\alpha|\leq1}\|f_{\alpha}\|_{H^{\frac{3}{2}-|\alpha|}(\gamma)}$
and  the constant $C$ only depends on $R$, $\theta$, $s$, $\gamma$,
$G$, $\|\psi\|_{C^2(\overline{G})}$, $\|\mu\|_{C^1(\overline{B_R})}$
and $\|\lambda\|_{C^1(\overline{B_R})}$.
\end{theorem}

$Proof$. By inverse trace theorem, there exists a $u^*\in H^2(G)$,
such that,\[\partial^{\alpha}u^*|_{\gamma}=f_{\alpha},\hspace{0.7cm}
|\alpha|\leq1\] and \begin{equation}\|u^*\|_{H^2(G)}\leq
C\zeta_0\end{equation} for some constant $C$ depending on $G$ and
$\gamma$. Hence\begin{equation}\|\widetilde{L}u^*\|_{L^2(G)}\leq
C\zeta_0.\end{equation}

We set $v=u-u^*$. By $(2.5)$, $v$ satisfies
\begin{equation}\partial^{\alpha}v|_{\gamma}=0,\hspace{0.7cm}
|\alpha|\leq1\end{equation} and
\[\widetilde{L}v=-\widetilde{L}u^*, \ \ {\rm{in}}\ G.\]

Set \[\psi(x)=R^2-|x|^2,\ \ {\rm{in}}\ G,\] and\[\varphi=e^{s\psi}-1.\]
Obviously, $\psi\in C^2(\overline{G})$ and
\begin{equation}\nabla\psi=-\left(
                             \begin{array}{c}
                               2x_1 \\
                               2x_2 \\
                               2x_3 \\
                             \end{array}
                           \right)
\neq0,\ \  {\rm{in}}\ G.\end{equation} We have \[\min_{x\in
\overline{G}}\varphi(x)=0,\]
\begin{equation}\varphi^*:=\max_{x\in
\overline{G}}\varphi(x)=e^{s(R^2-\theta^2)}-1>0.\end{equation}

Then by Theorem 2.1 and the fact of $C_c^{\infty}(G)$ being  dense
in $H_0^2(G)$, we know that there exist two positive constants $s_0$ and
$C$ such that for $s>s_0$, $\tau>\tau_0(s)$
\begin{equation}
\begin{array}{ll}
\displaystyle\tau^2s^4\int_{G}e^{2s\psi}e^{2\tau\varphi}|w|^2dx+&\displaystyle
s^2\int_{G}e^{2\tau\varphi}|\nabla
w|^2dx\\&+\displaystyle\frac{1}{\tau^2}\int_{G}e^{-2s\psi}e^{2\tau\varphi}|\nabla^2
w|^2dx\leq C\int_{G}e^{2\tau\varphi}|Lw|^2dx
\end{array}
\end{equation}
for any $w\in H^2_0(G)$.

We define the family $\{\omega(\delta)\}_{0<\delta<\varphi^*}$ of
subsets of $G$ by
\begin{equation}
\omega(\delta)=\{x\in G:\ \varphi(x)>\delta\}.
\end{equation}
Then the family satisfies
\[\emptyset=\omega(\varphi^*)\subset\omega(\delta')\subset\omega(\delta)\subset
\omega(0)=G\] for $0<\delta<\delta'<\varphi^*.$ Moreover,  it is easy to see
that $\omega(\delta)$ is a sub-domain of $G$ and we have
\[\partial \omega(\delta)\supset \gamma,\]
for each $ 0<\delta<\varphi^*$. Let
\begin{equation*}\eta>2.\end{equation*} Note \begin{equation*}
0<\mu=\frac{(1-\frac{1}{\eta-1})\varphi^*}{2}<\frac{\varphi^*}{2},\end{equation*}
and
 \begin{equation*}
\omega(\mu)\subset\subset\omega(\frac{\mu}{2}).\end{equation*}

Let $\chi\in C_c^{\infty}(B_R)$ be a cut off function
satisfying{\renewcommand{\baselinestretch}{2} \large\normalsize
\begin{equation}
\begin{array}{ll}
\displaystyle0\leq\chi\leq1,\ \ &  {\rm in}\ B_R,\\ \displaystyle
\chi=1\, &\displaystyle{\rm in}\ \omega(\mu)\cup \overline{B}_{\theta},\\
\chi=0\, &\displaystyle {\rm in}\ G\setminus \omega(\frac{\mu}{2}),\\
|\partial^\alpha\chi|<C_1,\ &{\rm in}\ B_R,\ |\alpha|\le1,
\end{array}
\end{equation}}
where $C_1$ is a constant depending on $R$ and the radius of the
support of $\chi$. Then by $(2.9)$ and $(2.14)$, we have $\chi v\in
H^2_0(G)$. Putting $\chi v$ into $(2.12)$ yields
\begin{equation*}
\begin{array}{ll}
\displaystyle\tau^2s^4\int_{G}e^{2s\psi}e^{2\tau\varphi}|\chi
v|^2dx+&\displaystyle s^2\int_{G}e^{2\tau\varphi}|\nabla \chi
v|^2dx\\&+\displaystyle\frac{1}{\tau^2}\int_{G}e^{-2s\psi}e^{2\tau\varphi}|\nabla^2
\chi v|^2dx\leq C\int_{G}e^{2\tau\varphi}|L(\chi v)|^2dx,
\end{array}
\end{equation*}
where $C$ depends on $\|\psi\|_{C^2(G)}$,  $\|\lambda\|_{C^1(G)}$
and $\|\mu\|_{C^1(G)}$. Ignoring  the second order term of the
left, one has
\begin{equation}
\displaystyle\tau^2s^4\int_{G}e^{2s\psi}e^{2\tau\varphi}|\chi v|^2dx+
s^2\int_{G}e^{2\tau\varphi}|\nabla \chi v|^2dx\leq
C\int_{G}e^{2\tau\varphi}|L(\chi v)|^2dx.
\end{equation}
Note that by $(2.3)$,{\renewcommand{\baselinestretch}{2}
\large\normalsize
\begin{eqnarray}
L(\chi v)&=&\chi Lv+[L,\chi]v\nonumber\\
&=&\chi{\Big{(}}\widetilde{L}v -{\big{(}}\nabla v+(\nabla
v)^T{\big{)}}\nabla\mu-{\rm{div}}\ v\nabla
\lambda{\Big{)}}+[L,\chi]v\nonumber\\
&=&\chi{\Big{(}}-\widetilde{L}u^* -{\big{(}}\nabla v+(\nabla
v)^T{\big{)}}\nabla\mu-{\rm{div}}\
v\nabla\lambda{\Big{)}}+[L,\chi]v,
\end{eqnarray}}
where the communicator $[L,\chi]v=L(\chi v)-\chi Lv$ is a system of
first order operators whose coefficients vanish on
${\big{(}}\omega(\mu)\cup \overline{B}_{\theta}{\big{)}}\cup
{\big{(}}G\setminus \omega(\frac{\mu}{2}){\big{)}} $. And the
coefficients of $[L,\chi]$ are bounded by a constant depending only
on $\lambda,\mu$ and $C_1$.

Then by $(2.15)$, $(2.16)$ and the triangle inequality,  we have
{\renewcommand{\baselinestretch}{2.5} \large\normalsize
\begin{equation*}\begin{array}{lll}
&&\displaystyle\tau^2s^4\int_{G}e^{2s\psi}e^{2\tau\varphi}|\chi v|^2dx+
s^2\int_{G}e^{2\tau\varphi}|\chi \nabla
v|^2dx-C_2s^2\int_{\omega(\frac{\mu}{2})\setminus\omega(\mu)}e^{2\tau\varphi}|v|^2dx\\&&\leq
\displaystyle C{\Big{(}}\int_{G}e^{2\tau\varphi}(|\chi
\widetilde{L}u^*|^2+|\chi\nabla
v|^2)dx+\int_{\omega(\frac{\mu}{2})\setminus\omega(\mu)}e^{2\tau\varphi}(|\nabla
v|^2+|v|^2)dx{\Big{)}}.\end{array}
\end{equation*}}
Choose $s$ large enough so that $C\displaystyle\int_G
e^{2\tau\varphi}|\chi\nabla v|^2dx$ can be absorbed into the left
side and move the last term of the left side to the right. Noting
that $\min\limits_{\overline{G}}\psi(x)=0$, we have
$e^{2s\psi}\geq1$ and
consequently{\renewcommand{\baselinestretch}{2.5} \large\normalsize
\begin{equation*}\begin{array}{lll}
&&\displaystyle\tau^2\int_{G}e^{2\tau\varphi}|\chi v|^2dx+
\int_{G}e^{2\tau\varphi}|\chi \nabla v|^2dx\\&&\leq \displaystyle
C{\Big{(}}\int_{G}e^{2\tau\varphi}|\chi
\widetilde{L}u^*|^2dx+\int_{\omega(\frac{\mu}{2})\setminus\omega(\mu)}e^{2\tau\varphi}(|\nabla
v|^2+|v|^2)dx{\Big{)}}.\end{array}
\end{equation*}}
Noting $\omega(\frac{\varphi^*}{2})\subset\omega(\mu)$ and by the
definition of $\chi$, we have{\renewcommand{\baselinestretch}{2.5}
\large\normalsize
\begin{equation}\begin{array}{lll}
&&\displaystyle\tau^2\int_{\omega(\varphi^*/2)}e^{2\tau\varphi}| v|^2dx+
\int_{\omega(\varphi^*/2)}e^{2\tau\varphi}| \nabla v|^2dx\\&&\leq
\displaystyle C{\Big{(}}\int_{G}e^{2\tau\varphi}|
\widetilde{L}u^*|^2dx+\int_{\omega(\frac{\mu}{2})\setminus\omega(\mu)}e^{2\tau\varphi}(|\nabla
v|^2+|v|^2)dx{\Big{)}}.\end{array}
\end{equation}}
Furthermore, the definition of $\{\omega(\delta)\}$ shows
that{\renewcommand{\baselinestretch}{2} \large\normalsize
\begin{equation}
\begin{array}{ll}
\displaystyle\varphi\geq\frac{\varphi^*}{2}\ \  & \displaystyle\rm{on}\ \omega(\frac{\varphi^*}{2}),\\
\varphi\leq\varphi^*\ &\rm{in}\ G,\\
\displaystyle\frac{\mu}{2}\leq\varphi\leq\mu\ &\displaystyle\rm{on}\  \omega(\frac{\mu}{2})\setminus\omega(\mu).\\
\end{array}
\end{equation}}
We obtain from $(2.17)$ and $(2.18)$ that
\begin{equation}\begin{array}{lll}
&&\displaystyle\tau^2e^{\tau\varphi^*}\int_{\omega(\frac{\varphi^*}{2})}|
v|^2dx+ e^{\tau\varphi^*}\int_{\omega(\frac{\varphi^*}{2})}| \nabla
v|^2dx\\&&\leq \displaystyle C{\Big{(}}e^{2\tau\varphi^*}\int_{G}|
\widetilde{L}u^*|^2dx+e^{2\tau\mu}\int_{\omega(\frac{\mu}{2})\setminus\omega(\mu)}(|\nabla
v|^2+|v|^2)dx{\Big{)}}.\end{array}
\end{equation}
Dividing both side of $(2.19)$ by $e^{\tau\varphi^*}$ and noting that
$v=u-u^*$, we have {\renewcommand{\baselinestretch}{2.5}
\large\normalsize
\begin{equation}\begin{array}{lll}
\displaystyle\int_{\omega(\frac{\varphi^*}{2})}(| v|^2+|\nabla
v|^2)dx&\leq& \displaystyle C{\Big{(}}e^{\tau\varphi^*}\int_{G}|
\widetilde{L}u^*|^2dx+e^{-\frac{\tau\varphi^*}{\eta-1}}\int_{\omega(\frac{\mu}{2})\setminus\omega(\mu)}(|\nabla
v|^2+|v|^2)dx{\Big{)}}\\&\leq&\displaystyle
C{\Big{(}}e^{\tau\varphi^*}\int_{G}|
\widetilde{L}u^*|^2dx+e^{-\frac{\tau\varphi^*}{\eta-1}}\int_{G}(|\nabla
v|^2+|v|^2)dx{\Big{)}}\\&\leq&\displaystyle
C{\Big{(}}e^{\tau\varphi^*}\int_{G}|
\widetilde{L}u^*|^2dx+e^{-\frac{\tau\varphi^*}{\eta-1}}(\|u\|_{H^1(G)}^2+\|u^*\|_{H^1(G)}^2){\Big{)}}\\
&\leq&C{\Big{(}}e^{\tau\varphi^*}\zeta_0^2+e^{-\tau\frac{\varphi^*}{\eta-1}}(M_0^2+\zeta_0^2){\Big{)}}.\end{array}
\end{equation}}

$Case$ 1: $\displaystyle M_0>
\zeta_0\exp\{\frac{3\tau_0(s)\varphi^*}{2(1-\frac{1}{\eta})}\}$.
Then $(2.20)$ becomes{\renewcommand{\baselinestretch}{2}
\large\normalsize
\begin{equation}\begin{array}{lll}
\displaystyle\int_{\omega(\frac{\varphi^*}{2})}(| v|^2+|\nabla
v|^2)dx&\leq& C{\Big{(}}e^{\tau\varphi^*}\zeta_0^2+e^{-\tau\frac{\varphi^*}{2}}(1+\exp\{-\frac{3\tau_0(s)\varphi^*}{2(1-\frac{1}{\eta})}\})M_0^2{\Big{)}}\displaystyle \\
&\leq&C(1+\exp\{-\frac{3\tau_0(s)\varphi^*}{2(1-\frac{1}{\eta})}\}){\Big{(}}e^{\tau\varphi^*}\zeta_0^2+e^{-\tau\frac{\varphi^*}{2}}M_0^2{\Big{)}}.\end{array}
\end{equation}}
Noting $\varphi^*>0$, we can put
\[\tau=2(1-\frac{1}{\eta})\frac{1}{\varphi^*}\ln\frac{M_0}{\zeta_0}.\]
It is easy to check that  $\tau>\tau_0(s)$, then the Carleman
estimate applies.  Choosing $$\omega=\omega(\frac{\varphi^*}{2})$$
and putting $\tau$ into $(2.21)$,
\[\|v\|_{H^1(\omega)}\leq
CM_0^{1-\frac{1}{\eta}}\zeta_0^{\frac{1}{\eta}}.\] Hence
\[\|u\|_{H^1(\omega)}\leq
CM_0^{1-\frac{1}{\eta}}\zeta_0^{\frac{1}{\eta}}
+\|u^*\|_{H^1(\omega)}\leq
CM_0^{1-\frac{1}{\eta}}\zeta_0^{\frac{1}{\eta}}.
\]
where $C$ only depends on $R$, $\theta$, $s$, $\gamma$, $G$,
$\|\psi\|_{C^2(\overline{G})}$ $\|\mu\|_{C^1(\overline{B_R})}$ and
$\|\lambda\|_{C^1(\overline{B_R})}$.

 $Case$ 2: $\displaystyle M_0\le
\zeta_0\exp\{\frac{3\tau_0(s)\varphi^*}{2(1-\frac{1}{\eta})}\}$. Then
we trivially have
\[\|u\|_{H^1(G)}=M_0=M_0^{1-\frac{1}{\eta}}M_0^{\frac{1}{\eta}}\leq
\exp\{\frac{3\tau_0(s)\varphi^*}{2(\eta-1)}\}M_0^{1-\frac{1}{\eta}}\zeta_0^{\frac{1}{\eta}},\]

This ends  the proof of Theorem 2.2.\hfill $\square$
\\

In particular the theorem shows the (local) uniqueness of the
solution to the Cauchy problem for $\widetilde{L}u=0$.

{\section*{{\normalsize \bf 3. Three spheres inequalities}}}
\setcounter{section}{3}\setcounter{equation}{0}\setcounter{lemma}{0}\setcounter{proposition}{0}\setcounter{theorem}{0}\setcounter{remark}{0}

We turn now to prove the three spheres inequality for the Lam\'{e}
system of elasticity. To begin, we recall a rather well known
interior estimate for elliptic systems (see \cite{5}, Ch. 8, Th. 2.2
for example).
\\
\begin{theorem}Assume $\lambda,\mu\in C^1(\Omega)$. Let $u\in
H^2(\Omega)$ be a solution of $(1.1)$. Then for any $B_r\subset B_R$
with $0<r<R$, there exists a constant $C$ depending only on
$\lambda$, $\mu$, $R$ and $r$ such that \[\|u\|_{H^2(B_r)}\leq
C\|u\|_{L^2(B_R)}.\]
\end{theorem}
 For proving a quantitative estimate of unique continuation in Section 4, we have to show how $C$ depends on $r$. Checking the proof of the above theorem in
 \cite{5}, we easily know that
there exists a constant $C_0$ independent of $R$ and $r$ such that
the following estimate holds.
\[\|u\|_{H^2(B_r)}\leq\frac{C_0}{(R-r)^2}
\|u\|_{L^2(B_R)}.\]

Inspired by \cite{12} and \cite{15}, we give a proof of the three
spheres inequality (Theorem 1.1) as
follows.\\
{\bf{Proof of Theorem 1.1}}. From the proof of Theorem 2.2, we can
see that there exists a constant $\theta_1$ satisfying
$\theta<\theta_1<R$ such that{\renewcommand{\baselinestretch}{2}
\large\normalsize
\begin{eqnarray}\|u\|_{H^1(B_{\theta_1}\setminus B_{\theta})}&\leq&
C{\Big{(}}\sum_{|\alpha|\le1}\|u\|_{H^{\frac{3}{2}-|\alpha|}(\partial
B_{\theta})}{\Big{)}}^{\epsilon}\|u\|_{H^2(B_{R})}^{1-\epsilon}\nonumber\\
&\le&C\|u\|_{H^2(B_\theta)}^\epsilon\|u\|_{H^2(B_{R})}^{1-\epsilon},\end{eqnarray}}
where $\epsilon\in(0,1)$. Furthermore,
\[\|u\|_{H^1(B_{\theta_1})}-\|u\|_{H^1( B_{\theta})}\le \|u\|_{H^1(B_{\theta_1}\setminus B_{\theta})}\le C\|u\|_{H^2(B_\theta)}^\epsilon\|u\|_{H^2(B_{R})}^{1-\epsilon}.\]
Then we have
\begin{equation}\|u\|_{L^2(B_{\theta_1})}\le \|u\|_{H^1(B_{\theta_1})}\le
C\|u\|_{H^2(B_{\theta})}^{\epsilon}\|u\|_{H^2(B_{R})}^{1-\epsilon}.\end{equation}

Setting $\theta_2=\frac{\theta+\theta_1}{2}$ and using Theorem 3.1,
we obtain
\begin{equation*}\|u\|_{H^2(B_{\theta_2})}\leq \frac{C}{(\theta_1-\theta)^2}
\|u\|_{L^2(B_{\theta_1})}.\end{equation*} Combing $(3.2)$, we have
\begin{equation}\|u\|_{H^2(B_{\theta_2})}\leq
C\|u\|_{H^2(B_{\theta})}^{\epsilon}\|u\|_{H^2(B_{R})}^{1-\epsilon}.\end{equation}
Set \[\displaystyle R_0=\frac{R+R_2}{2}\in(R_2,R),\] \[\displaystyle
\theta:=\frac{R_2}{R_0}R,\ \displaystyle
a:=\frac{\theta_2}{\theta}>1.\]

Claim that\begin{equation} \|u\|_{H^2(B_{ra})}\leq
C{\Big{(}}\frac{\theta^2}{r^2}\|u\|_{H^2(B_{R_0})}^{1-\epsilon}{\Big{)}}\|u\|_{H^2(B_{r})}^{\epsilon},
\ \ {\rm{for \ all}}\ 0<r<R_2.
\end{equation}
where $C$ and $\epsilon$ depend only on $R$, $R_2$, $\theta$ and the
Lam\'{e} moduli $\lambda$, $\mu$. Setting $\displaystyle
t:=\frac{r}{\theta }<1$, we have\begin{equation} r=\theta
t<ra=\theta_2t<R_0<R.\end{equation}

Let $\widetilde{u}(y):=u(ty)$,
$\widetilde{\lambda}(y):=\lambda(ty)$,
$\widetilde{\mu}(y):=\mu(ty)$. Then $\widetilde{u}$ satisfies

\begin{equation}
\rm{div}{\big{(}}\widetilde{\mu}(\nabla \widetilde{u}+(\nabla
\widetilde{u})^\top){\big{)}}+\nabla{\big{(}}\widetilde{\lambda}\rm{div}\
\widetilde{u}{\big{)}}=0,
\end{equation}
where $\widetilde{\lambda}, \widetilde{\mu}\in C^1(\Omega)$ satisfy
strong ellipticity conditions. Repeating the same argument, we have
\begin{equation}\|\widetilde{u}\|_{H^2(B_{\theta_2})}\leq
C\|\widetilde{u}\|_{H^2(B_{\theta})}^{\epsilon}\|\widetilde{u}\|_{H^2(B_{R})}^{1-\epsilon}.\end{equation}
The constant $C$ appearing in $(3.7)$ is independent of $t$.

Indeed, $\widetilde{u}$ satisfies $(3.6)$ in $B_{\frac{R}{t}}$.
Noting that $t<1$, we deduce $\widetilde{u}$ also satisfies $(3.6)$
in $B_{R}$. Then applying Theorem 2.2, we know the constant $C$ in
$(3.7)$ depends on $\theta$, $R_2$, $R$,
$\|\widetilde{\lambda}\|_{C^1(\overline{B_R})}$,
$\|\widetilde{\mu}\|_{C^1(\overline{B_R})}$, and
$\|\psi\|_{C^2(\overline{B_R\backslash B_{\theta}})}$, while the
coefficients satisfy $\widetilde{\mu}\geq\alpha_0>0,$
$2\widetilde{\mu}+\widetilde{\lambda}\geq \beta_0>0$ and
\[\sup_{B_R}\{|\partial^\alpha\widetilde{\mu}|,|\partial^\alpha\widetilde{\lambda}|\}\leq
\sup_{B_{R}}\{|\partial^\alpha\mu|,|\partial^\alpha\lambda|\}.\]

By $t<1$ and the change of variable $ty=x$, we have
\[\|u\|_{H^2(B_{ra})}\leq
C\frac{1}{t^2}\|u\|_{H^2(B_{r})}^{\epsilon}\|u\|_{H^2(B_{R_0})}^{1-\epsilon},\]
which gives $(3.4)$.

Choose $\displaystyle r=\frac{R_1}{2}$. Then there exists a unique
positive integer $N$ such that \begin{equation}ra^{N-1}<R_2\leq
ra^N.\end{equation} Since $\displaystyle
ra^N<aR_2=\frac{\theta_2}{\theta}R_2<R_0$, we have $ra^N<R$ and
\begin{equation}\|u\|_{H^2(B_{ra^k})}\leq
E_k\|u\|_{H^2(B_{ra^{k-1}})}^{\epsilon},\end{equation} for all
$k=1,2,\cdots,N$, where we
set\[E_k:=C(\frac{\theta}{ra^{k-1}})^2\|u\|_{H^2(B_{R_0})}^{1-\epsilon}.\]
Since $a>1$, we have $E_k<E_1$ for $k=2,\cdots,N$. Then repeated use
of $(3.6)$ shows that {\renewcommand{\baselinestretch}{2}
\large\normalsize
\begin{eqnarray}
\|u\|_{H^2(B_{R_2})}&\leq&\displaystyle \|u\|_{H^2(B_{ra^N})}\nonumber\\
&\leq&\displaystyle E_N\|u\|_{H^2(B_{ra^{N-1}})}^{\epsilon}\nonumber\\
&\leq&\displaystyle E_1{\Big{(}}E_{N-1}\|u\|_{H^2(B_{ra^{N-2}})}^{\tau_1}{\Big{)}}^{\epsilon}\nonumber\\
&\leq&\displaystyle
E_1^{\frac{1-\epsilon^N}{1-\epsilon}}\|u\|_{H^2(B_{r})}^{\epsilon^N}\nonumber
\end{eqnarray}}

Setting $\sigma=\epsilon^N<1$, we have
\begin{equation*}
\|u\|_{H^2(B_{R_2})}\leq
C\|u\|_{H^2(B_r)}^{\sigma}\|u\|_{H^2(B_{R_0})}^{1-\sigma}.
\end{equation*}
By Theorem 3.1, we obtain
\[\|u\|_{H^2(B_{R_2})}\leq
C\frac{1}{R_1^2}\|u\|_{L^2(B_{R_1})}^{\sigma}\|u\|_{L^2(B_{R})}^{1-\sigma}.\]
This implies \begin{equation}\|u\|_{L^2(B_{R_2})}\leq
C\frac{1}{R_1^2}\|u\|_{L^2(B_{R_1})}^{\sigma}\|u\|_{L^2(B_{R})}^{1-\sigma},\end{equation}
where $C$ and $\sigma$ depend on $\frac{R_1}{R_2}$, $\frac{R_2}{R}$
and coefficients $\lambda$, $\mu$.
We thus  complete
the proof.\hfill $\square$ \ \\

\begin{remark}
Furthermore, we can easily obtain, from $(3.4)$ and $(3.10)$ that
\begin{equation} C\leq\frac{\widetilde{C}}{R_1^4},\end{equation}
where $\widetilde{C}$ depends on $R_2$, $R$,
$\|\lambda\|_{C^1(\overline{B_R})}$,
$\|\mu\|_{C^1(\overline{B_R})}$,
$\|\psi\|_{C^2(\overline{B_R\backslash B_{\theta}})}$. By the
definition of $\psi$, we know that
$\|\psi\|_{C^2(\overline{B_R\backslash B_{\theta}})}$ depends on $R$
and $\theta(=2\frac{R_2R}{R_2+R})$.
\end{remark}

We  are now at a position to discuss SUCP.  The strong unique
continuation is close related with the three spheres  inequality
(see \cite{4} for the case of scalar parabolic equations). However,
the three spheres inequality obtained in this paper may not be used
to prove SUCP of Lam\'{e} systems, since the constant $C$ and
$\sigma$ appeared in the right hand of $(1.2)$ both depend on $R_1$.
Fortunately, we can get a weak sense of SUCP.

{\section*{{\normalsize \bf 4.  Unique continuation}}}
\setcounter{section}{4}\setcounter{equation}{0}\setcounter{lemma}{0}\setcounter{proposition}{0}\setcounter{theorem}{0}

\ \\{\bf{Proof of Theorem 1.2}}. {\rm Part (i)}.  Without loss of
generality we may assume that $x_0=0$, i.e.,
\begin{equation*}\int_{B_r}|u|^2dx=O(e^{-{r^{-\varepsilon}}}),\ \ {\rm{as}}\ r\rightarrow0.\end{equation*} We wish
to show that $u\equiv0$ in $\Omega$.

The following proof is based on the proof of Theorem 1.1 and Theorem
2.2, from which we utilize the notation and terminology.

By $(3.8)$, we know\[\frac{R_1}{2}a^{N-1}<R_2\leq
\frac{R_1}{2}a^N.\] Hence\begin{equation} (\ln
a)^{-1}\ln\frac{2R_2}{R_1}\leq N<(\ln
a)^{-1}\ln\frac{2R_2}{R_1}+1\end{equation} where $a>1$ and $N$ are
defined in the proof of Theorem 1.1.

$(1.4)$ implies,
\begin{equation} \int_{B_{R_1}}|u|^2dx\leq
Ce^{-{R_1^{-\varepsilon}}},\ \ 0<R_1<1.\end{equation} In order to
prove $u\equiv 0$, we need to find a proper $\epsilon$ which appears
in $(2.6)$. Let
\begin{equation}\eta:=\exp{\Big{\{}}\frac{1}{N}+1{\Big{\}}}, \end{equation}
 Note
that $\eta>2$.

Repeating the same discussion as the proof of Theorem 2.2, we have
\begin{equation*}\|u\|_{H^1(\omega)}\leq
CM_0^{1-\frac{1}{\eta}}\zeta_0^{\frac{1}{\eta}},\end{equation*}
where $C$ is independent of $R_1$. By the proof of Theorem 1.1, we
know
\begin{equation*}\|u\|_{L^2(B_{R_2})}\leq
C\|u\|_{L^2(B_{R_1})}^\sigma\|u\|_{L^2(B_{R})}^{1-\sigma},\end{equation*}
where $\displaystyle\sigma=\frac{1}{\eta^N}$. Noting $(4.1)$ and
$(4.3)$, we have\[\sigma=\frac{1}{\eta^{N}}=e^{-(1+N)}\ge
e^{-2}{\Big{(}}\frac{R_1}{2R_2}{\Big{)}}^{(\ln a)^{-1}}.\] By
$(4.1)$ and $(4.2)$, we have{\renewcommand{\baselinestretch}{2.5}
\large\normalsize
\begin{eqnarray}\int_{B_{R_2}}|u|^2dx&\leq&\displaystyle
C{\Big{(}}\int_{B_{R}}|u|^2dx{\Big{)}}^{1-\sigma}{\Big{(}}\int_{B_{R_1}}|u|^2dx{\Big{)}}^{\sigma}\nonumber\\
&\leq& \displaystyle C e^{-{\sigma R_1^{-\varepsilon}}}{\Big{(}}\int_{B_{R}}|u|^2dx{\Big{)}}^{1-\sigma}\nonumber\\
&\leq& \displaystyle \frac{\widetilde{C}}{R_1^4}
\exp\{e^{-2}\frac{-1}{R_1^{\varepsilon-(\ln
a)^{-1}}}\}{\Big{(}}\int_{B_{R}}|u|^2dx{\Big{)}}^{1-\sigma}.
\end{eqnarray}}

~~~\\Claim: $(\ln a)^{-1}<\varepsilon.$

Indeed, from the proof of Theorem 1.1 and Theorem 2.2, we know
$$a=\frac{\theta_2}{\theta}=\frac{1}{2}(1+\frac{\theta_1}{\theta})$$ where $\theta\le\theta_1\le R$ such that $B_{\theta_1}=\omega=\{x\in G:\ \varphi(x)>\frac{\varphi^*}{2}\}.$
Noting the definition of $\varphi$:
\[\varphi(x)=e^{s(R^2-|x|^2)}-1,\]
we know $\theta_1$ does not depend on $\theta$.

The claim follows as long as  we let $\theta$ be small enough.\\

Then we pass to the limit in $(4.4)$ as $R_1\rightarrow0$,
\begin{equation*} \|u\|_{L^2(B_{R_2})}\leq 0.\end{equation*}
This implies \[u\equiv0,\ \ {\rm{in}}\ B_{R_2}.\]  {\rm Part (i)}
follows
by standard arguments.\\

\rm{Part (ii)}. Since $\lambda,\mu\in C^2(\Omega)$, $\partial_s u$
$(s=1,\cdots,n)$ also satisfies the Lam\'{e} system $(1.1)$, by the
same argument we have
\[\int_{B_{R_2}}|\partial_s u|^2dx\le
C{\Big{(}}\int_{B_{R}}|\partial_su|^2dx{\Big{)}}^{1-\sigma}{\Big{(}}\int_{B_{R_1}}|\partial_su|^2dx{\Big{)}}^{\sigma}.\]
This implies $\partial_s u=0,$ $(s=1,\cdots,n)$. Then Theorem 1.2 follows.

\hfill$\square$
\\

{\bf{Acknowledgment.}} The author wants to show his great gratitude
to Professor Hongwei Lou, Professor Xu Zhang and Professor Ching-Lung Lin for many helpful
guidance and suggestions.
\\\\


\footnotesize\begin{thebibliography}{99}

\bibitem{1}G. Alessandrini, A. Morassi, Strong unique continuation for the
Lam\'{e}system of elasticity, Comm. in. PDE, 26 (2001), pp.
1787-1810.
\bibitem{2}G. Alessandrini, A. Morassi, E. Rosset, Detecting an inclusion in an elastic body by boundary measurements, SIAM Review, 46 (2004), pp. 477-498.
\bibitem{3}  D. D. Ang, M. Ikehata , D. D. Trong, M. Yamamoto,  Unique continuation for a stationary
isotropic Lam\'{e}system with varaiable coefficients, Comm. in. PDE,
23 (1998), pp. 371-385.
\bibitem{4} B. Canuto, E. Rosset, S. Vessella, Quantitative estimates of unique continuation
for parabolic equations and inverse initial-boundary value problems
with unknown boundaries, Trans. Amer. Math. Soc., 354 (2002), pp.
491-535.
\bibitem{5}Y. -Z. Chen, L. -C. Wu, Second Order Elliptic Equations and Elliptic
Systems, Translations of Mathematical Monographs 174, AMS,
Providence, RI, 1998.

\bibitem{6}B. Dehman, L. Robbiano, La propri\'{e}t\'{e} du prolungement unique
pour un syst`eme elliptique: le syst\`{e}me de Lam\'{e}, J. Math.
Pures Appl., 72 (1993), pp. 475-492.
\bibitem{7}M. Eller, Carleman estimates for some elliptic systems,
Journal of Physics: Conference Series, 124 (2008) 012023.
\bibitem{8}M. Eller, V. Isakov, G. Nakamura,  D. Tataru, Uniqueness and
stability in the Cauchy problem for Maxwell and elasticity systems,
in Nonlinear Partial Differential Equations and Their Applications,
College de France Seminar, Vol. 14, Stud. Math. Appl. 31,
North¨CHolland, Amsterdam, 2002, pp. 329-349.
\bibitem{9}L. Escauriaza, Unique continuation for the system of
elasticity in the plan, Proc. Amer. Math. Soc., Vol. 134, 7 (2005),
pp. 2015-2018.
\bibitem{10}
N. Garofalo, F. H. Lin, Monotonicity properties of variational
integrals, Ap weights and unique continuation, Indiana Univ. Math.
J., 35 (1986), pp. 245-268.

\bibitem{11} N. Garofalo, F. H. Lin, Unique continuation for elliptic
operators: A geometricvariational approach, Comm. Pure Appl. Math.,
40 (1987), pp. 347-366.

\bibitem{12}N. Higashimori, A conditional stability estimate for identifying a cavity by an elastostatic measurement, J. Inverse
Ill-Posed Probl., in press.

\bibitem{13}I. Kukavika, Quantitative uniqueness for second-order elliptic
operators, Duke Math. J., 91 (1998), pp. 225-240.
\bibitem{14}E. M. Landis, A three-sphere theorem, Dokl. Akad. Nauk SSSR 148 (1963) 277¨C279, translated in Sov. Math. 4
(1963), pp. 76-78.
\bibitem{15}C. -L. Lin, G. Nakamura, J. -N. Wang, Three spheres
inequlities for a two-dimensional elliptic system and its
application, J. Differential Equations, 232 (2007), pp. 329-351.
\bibitem{16}C. -L. Lin, J. -N. Wang, Strong unique continuation for the Lam\'{e}
 system with Lipschitz coefficients, Math. Ann., 331 (2005), pp. 611-629.
\bibitem{17}L. E. Payne,  H. F. Weinberger,  New bounds for solutions
of second order elliptic partial differential equations, Pacific J.
Math., 8 (1958), pp. 551-573.

\bibitem{18}N. Weck, Au{\ss}nraumaufgaben in der Theorie station\"{a}rer Schwingungen inhomogener
elasticher K\"{o}rper, Math. Z., 111 (1969), pp. 387-398.
\bibitem{19}N. Weck, Unique continuation for systems with Lam\'{e} principal part, Math. Methods
Appl. Sci., 24 (2001), pp. 595-605.



\end{thebibliography}
\end{document}